\documentclass[a4paper,11pt]{amsart}
\usepackage{amsmath,amsthm,amssymb}
\usepackage[all]{xy}
\DeclareMathOperator{\sgn}{sgn}
\DeclareMathOperator{\idem}{idem}
\DeclareMathOperator{\pol}{pol}
\DeclareMathOperator{\pert}{pert}
\DeclareMathOperator{\Rea}{Re}
\DeclareMathOperator{\Skvol}{{\mathfrak{Skvol}}}
\DeclareMathOperator{\Invol}{{\mathfrak{Invol}}}
\DeclareMathOperator{\Idem}{{\mathfrak{Idem}}}
\DeclareMathOperator{\spec}{{\mathfrak{Spec}}}

\newcommand{\qand}{\qquad\text{and}\qquad}
\theoremstyle{definition}
\newtheorem{point}{}[section]
\newtheorem{defin}[point]{Definition}

\newtheorem{remark}[point]{Remark}
\newtheorem{objective}[point]{Objective}
\theoremstyle{plain}
\newtheorem{prop}[point]{Proposition}
\newtheorem{lemma}[point]{Lemma}

\newcommand{\marginextend}[1]{ \addtolength{\oddsidemargin}{-#1}  \addtolength{\evensidemargin}{-#1}
  \addtolength{\textwidth}{#1}\addtolength{\textwidth}{#1}}
\newcommand{\updownextend}[1]{ \addtolength{\topmargin}{-#1}  \addtolength{\textheight}{#1}
\addtolength{\textheight}{#1}}
\marginextend{1.3cm}
\updownextend{0.35cm}
\allowdisplaybreaks[1]
\title{Spectral calculations in rings}
\author{Gyula Lakos}
\address{Department of Geometry, E\"otv\"os University, P\'azm\'any P\'eter s.~1/C,  Budapest, H--1117, Hungary}
\email{lakos@cs.elte.hu}
\keywords{Sign operation, square root operation, formal transformation kernels}
\subjclass[2000]{Primary: 46H30, Secondary:  13J99.}
\begin{document}
\begin{abstract}
We examine the validity of certain spectral integral formulas in topological rings.
We consider the sign and square-root functions in polymetric rings containing $\frac12$.
It turns out that formal analogues of classical transformation kernels
and the resolvent identity can be used to understand the situation.
In the lack of $\frac12$, the functions
$\frac12-\frac12\sgn\left(\frac12-z\right)$ and $\frac12-\sqrt{\frac14-z}$
can be generalized, respectively.
\end{abstract}
\maketitle
\setcounter{section}{-1}
\section{Introduction}
Spectral integrals like
\begin{equation}
\sgn Q=\int_{\{z\in\mathbb C\,:\, |z|=1\}}
\frac{\frac{1-z}2+\frac{1+z}2Q}{\frac{1+z}2+\frac{1-z}2Q}\,\frac{|dz|}{2\pi}
\label{e1}\end{equation}
and
\begin{equation}\sqrt S=\int_{\{z\in\mathbb C\,:\, |z|=1\}}
\frac{S}{\frac1z\left(\left(\frac{1+z}{2}\right)^2-
\left(\frac{1-z}{2}\right)^2S\right)}\,\frac{|dz|}{2\pi}
\label{e2}\end{equation}
are often useful.
They extend the complex functions $\sgn Q$, which is the sign function of (the real part of) $Q$,
and $\sqrt S$, which is the square root function cut along the negative real axis, respectively.
Definitions like above are justified if
they are supported by appropriate algebraic identities and spectral properties.
This is the situation in linear analysis, where the formulas above can be established
for elements with appropriate spectral properties
in great generality, even if the resolvent terms are not necessarily
continuous, cf. Haase \cite{H}, Mart\'\i nez Carracedo--Sanz Alix \cite{MS}.
We refer to this case as the ``analytic case''.

However, in the analytic case, if the resolvent terms are continuous, then they are also smooth, and
one can expand everything in terms of Fourier series, or rather Laurent series in  $z$.
One can naturally ask if similar computations can be done in more general rings,
in particular, in rings without a natural $\mathbb R$-action.
It is natural to  check these ideas for formal Laurent series on polymetric rings.
We refer to this case as the ``algebraic case''.

We call a  ring $\mathfrak A$  polymetric if
\begin{itemize}
\item[(a)] its topology is induced by a family of ``seminorms'' $p:\mathfrak A\rightarrow[0,+\infty)$ such that
$p(0)=0$, $p(-X)=p(X)$, $p(X+Y)\leq p(X)+p(Y)$,
\item[(b)]for each ``seminorm'' $p$ there exists a ``seminorm'' $\tilde p$ such that $p(XY)\leq\tilde p(X)\tilde p(Y)$;
\end{itemize}
i.~e., if it is a polymetric space whose multiplication is compatible with the topology.
Dealing with Laurent series in $z$, integration over the unit circle  becomes a formal process.
It is nothing else but detecting the coefficient of $z^0$.
On the other hand, for the multiplication of Laurent series, some sort of convergence control is required.
Primarily, we will be interested in Laurent series with rapidly decreasing coefficients.
A  sequence is rapidly decreasing if it is rapidly decreasing in each seminorm.
Furthermore, we consider only sequentially complete, Hausdorff polymetric rings.
We also assume that $\frac12\in\mathfrak A$.
Then \eqref{e1} and \eqref{e2} are meaningful.

Indeed, applied to elements with appropriate ``spectral'' properties, the expressions
$\sgn Q$ and $\sqrt S$ will have good properties justifying the notation.
This can be proved by an analysis of the coefficients.
The algebraic approach is particularly manageable in the case of \eqref{e1},
it is  essentially shown in Karoubi \cite{K}, by a direct analysis of coefficients,
that the expression  $\sgn Q$ yields an involution
compatible with factorization of affine loops.
Nevertheless, such computations are not necessarily very enlightening.
The objective of this paper is to prove our statements regarding the algebraic case and
to do this in a manner which brings the algebraic and analytic cases together, at least formally.
It turns out that the basic tool of the analytic case, the resolvent identity, works generally.
Another natural question is what happens in the lack of $\frac12$.
Then the sign and square root functions are not really appropriate.
Instead, we can generalize
\begin{equation}
\idem P=\int_{\{z\in\mathbb C\,:\, |z|=1\}}\frac{Pz}{1-P+Pz}\,\frac{|dz|}{2\pi}
\label{e3}\end{equation}
and
\begin{equation}
\sqrt[\mathcal F] T=\int_{\{z\in\mathbb C\,:\, |z|=1\}}
\frac{(1+z)T}{1+(z-2+z^{-1})T}\,\frac{|dz|}{2\pi},
\label{e4}\end{equation}
extending the functions
$\idem P=\frac12-\frac12\sgn\left(\frac12-P\right)$ and $\sqrt[\mathcal F] T=\frac12-\sqrt{\frac14-T}$,
repectively.

\section{Laurent series}
\begin{point}
If $\mathfrak A$ is a  polymetric  ring, then we may consider formal
Laurent series $a=\sum_{n\in\mathbb Z}a_nz^n$.
If $p$ is a seminorm on $\mathfrak A$ and
$\alpha:\mathbb Z\rightarrow\mathbb R_+$ is a non-negative
function, then we may define $p_\alpha(a)=\sum_{n\in\mathbb Z}\alpha(n)p(a_n)$.
Then we may consider the  polymetric spaces
(a) $\mathfrak A[z^{-1},z]^{\mathrm f}$ of essentially finite Laurent series,
(b) $\mathfrak A[z^{-1},z]^{\infty}$ of rapidly decreasing Laurent series,
(c) $\mathfrak A[z^{-1},z]^{\mathrm b}$ of ``summable'' Laurent series,
(d) $\mathfrak A[[z^{-1},z]]^{\mathrm b}$ of bounded Laurent series,
(e) $\mathfrak A[[z^{-1},z]]^{\infty}$ of polynomially growing Laurent series,
(f) $\mathfrak A[[z^{-1},z]]^{\mathrm f}$ of formal Laurent series,
as the spaces which contain series bounded for each $p_\alpha$  such that
(a) $\alpha$ is unrestricted,
(b) $\alpha$ is polynomially growing,
(c) $\alpha$ is bounded,
(d) $\alpha$ is summable,
(e) $\alpha$ is rapidly decreasing,
(f) $\alpha$ is vanishing except at finitely many places, respectively.
We have continuous inclusions
\[\mathfrak A[z^{-1},z]^{\mathrm f}\hookrightarrow \mathfrak A[z^{-1},z]^{\infty}\hookrightarrow
\mathfrak A[z^{-1},z]^{\mathrm b}\hookrightarrow\mathfrak A[[z^{-1},z]]^{\mathrm b}\hookrightarrow
\mathfrak A[[z^{-1},z]]^{\infty}\hookrightarrow \mathfrak A[[z^{-1},z]]^{\mathrm f}.\]

Of these spaces, $\mathfrak A[z^{-1},z]^{\mathrm f}$,
$\mathfrak A[z^{-1},z]^{\infty}$, $\mathfrak A[z^{-1},z]^{\mathrm b}$
will remain polymetric rings.
Indeed, $\widetilde{(p_\alpha)}$ can be chosen as $\tilde p_{\check\alpha}$, where
$\check\alpha(n)=1\vee\max_{-2|n|\leq m\leq 2|n|} |\alpha(m)|$.
We have compatible continuous module actions
$\mathfrak A[z^{-1},z]^{\bullet}\times \mathfrak A[[z^{-1},z]]^{\bullet}
\rightarrow\mathfrak A[[z^{-1},z]]^{\bullet}$.
Essentially the same applies to the spaces of power series
$\mathfrak A[z]^{\mathrm f}$, $\mathfrak A[z]^{\infty}$, $\mathfrak A[z]^{\mathrm b}$,
$\mathfrak A[[z]]^{\mathrm b}$, $\mathfrak A[[z]]^{\infty}$, $\mathfrak A[[z]]^{\mathrm f}$, except here even
$\mathfrak A[[z]]^{\infty}$, $\mathfrak A[[z]]^{\mathrm f}$ are  polymetric rings.
Indeed, for them,  $\widetilde{(p_\alpha)}$ can be chosen as $\tilde p_{\grave\alpha}$, where
$\grave\alpha(n)=\sqrt{\max_{n\leq m} |\alpha(m)|}$.
An element like $\frac{1+Q}2+\frac{1-Q}2z$  may be considered either as an element of
$\mathfrak A[z^{-1},z]^{\mathrm f}$, $\mathfrak A[z^{-1},z]^{\infty}$, or $\mathfrak A[z^{-1},z]^{\mathrm b}$, etc.
Practically, the difference is that
the larger the ring is the easier is to find a multiplicative inverse of the element given.
\end{point}
\begin{point}
If $a(z)=\sum_{n\in\mathbb Z}a_nz^n\in\mathfrak A[[z^{-1},z]]^{\mathrm f}$, then we define formally
\[\int a(z)\,\frac{|dz|}{2\pi}=a_0.\]

The Hilbert kernel (``up to multiplication by $\mathrm i$'') is defined as
\[\left[\frac{1+z}{1-z}\right]=\sum_{s\in\mathbb Z}(\sgn s)z^s\in\mathfrak A[[z^{-1},z]]^{\mathrm b}. \]
\end{point}
\begin{point}
Some further terminology is as follows.
For $a(z)=\sum_{n\in\mathbb Z}a_nz^n\in\mathfrak A[z^{-1},z]^{\mathrm b}$ we let
\[\lim_{z\nearrow1}a(z)=\sum_{n\in \mathbb Z}a_n.\]
Naturally, this notation also applies to power series.
\end{point}
In what follows, we let $\mathfrak A[z^{-1},z]$, $\mathfrak A[[z^{-1},z]]$
 $\mathfrak A[z]$, $\mathfrak A[[z]]$
denote
$\mathfrak A[z^{-1},z]^\infty$, $\mathfrak A[[z^{-1},z]]^\infty$,
$\mathfrak A[z]^\infty$, $\mathfrak A[[z]]^\infty$,
  respectively,  but
similar statements hold for ${}^{\mathrm f}$ and ${}^{\mathrm b}$, too.
\begin{prop} For  $a(t,z)\in\mathfrak A[t] \,[[z^{-1},z]]$
and $b(z)\in\mathfrak A[z^{-1},z]$,
\[\int\left(\lim_{t\nearrow1}a(t,z)\right)b(z)\,\frac{|dz|}{2\pi}
=\lim_{t\nearrow1}\int\left(a(t,z)\right)b(z)\,\frac{|dz|}{2\pi}.\]
\begin{proof}
This is just the generalized associativity of the rapidly decreasing (hence absolute convergent) sum
$\sum_{n\in\mathbb N,s\in\mathbb Z}a_{n,s}b_{-s}$.
\end{proof}
\end{prop}
\begin{point}
For the sake of brevity, we call the elements of $\mathfrak A[t]\,[[z^{-1},z]]$
as transformation kernels.
Practically, the convenient thing is to consider those elements of the \mbox{\textit{ring}}
$\mathfrak A[[t]]^{\mathrm f}[z^{-1},z]^{\mathrm b}$ which can  be thought to be transformation kernels.
(This is advantageous from computational viewpoint, because the product of $a(t,z)$ and
$b(z)\in\mathfrak A[z^{-1},z]$ formally yields the same element of
$\mathfrak A[[t]]^{\mathrm f}[[z^{-1},z]]^{\mathrm b}$
either we interpret $a(t,z)\in\mathfrak A[[t]]^{\mathrm f}[z^{-1},z]^{\mathrm b}$ or
$a(t,z)\in\mathfrak A[t]\,[[z^{-1},z]]$, but the first case is often easier to compute with.)

Such elements are the  Poisson kernel
\[\mathcal P(t,z)=\frac{1-t^2}{(1-tz)(1-tz^{-1})}=\sum_{s\in\mathbb Z} t^{|s|}z^s,\]
the Hilbert-Poisson kernel
\[\mathcal H(t,z)=\frac{t(z-z^{-1})}{(1-tz)(1-tz^{-1})}=\sum_{s\in\mathbb Z}(\sgn s)t^{|s|}z^s,\]
the $\frac12$-shifted odd Poisson kernel
\[\mathcal L(t,z)=\frac{(1-t)(1+z)}{(1-tz)(1-tz^{-1})}=\sum_{s\in\mathbb N}t^sz^{-s}+
\sum_{s\in\mathbb N}t^sz^{s+1},\]
and the variant regularization kernel
\[\tilde{\mathcal R}(t,z)=\frac{(1+t)t(1-z)(1-z^{-1})}{2(1-tz)(1-tz^{-1})}=
\sum_{s\in\mathbb N, s>0} \frac{t^{s+1}-t^s}{2}z^{-s}+t+\sum_{s\in\mathbb N, s>0} \frac{t^{s+1}-t^s}{2}z^s.\]
This latter one has the property
\[\lim_{t\nearrow1} \tilde{\mathcal R}(t,z)=1\]
(here we think of  $\tilde{\mathcal R}(t,z)$ as an element of $\mathfrak A[t]\,[[z^{-1},z]]$).
The ordinary regularization kernel
\[\mathcal R(t,z)=\frac{t(1-z)(1-z^{-1})}{(1-tz)(1-tz^{-1})}\]
is just an element of $\mathfrak A[[t]]^{\mathrm f}[z^{-1},z]^{\mathrm b}$, hence
it is not so convenient algebraically.
On the other hand, the use of variant regularization makes the variant Hilbert-Poisson kernel
\[\tilde{\mathcal H}(t,z)=\frac{1+t}2\mathcal H(t,z)=
\frac{(1+t)t(z-z^{-1})}{2(1-tz)(1-tz^{-1})}=\sum_{s\in\mathbb Z}(\sgn s)\frac{1+t}2t^{|s|}z^s\]
useful.\end{point}
It is not hard to see that the definitions and the proposition above can be formulated in the case
when we have many variables $z,w,\ldots$ instead of just $z$. For example, we may consider the Hilbert kernel
\[\left[\frac{z+w}{z-w}\right]=\left[\frac{1+wz^{-1}}{1-wz^{-1}}\right]=\sum_{s\in\mathbb Z}(\sgn s) w^sz^{-s}.\]

\section{Spectral classes}
\begin{point}
In order to save some space we use the short-hand notation
\begin{align}
\Lambda(a)&=\tfrac{1}{2}(1+a),\notag\\
\Lambda(a,b)&=\tfrac{1}{2}(1+a+b-ab),\notag\\
\Lambda(a,b,c)&=\tfrac{1}{4}(1+a+b+c-ab+ac-bc+abc),\notag\\
\Lambda(a,b,c,d)&=\tfrac{1}{4}(1+a+b+c+d-ab-bc-cd+ac+ad+bd\notag\\&
  \qquad+abc-acd-abd+bcd-abcd),\notag
\end{align}
etc., following the scheme
\[\Lambda(c_1,c_2,c_3,\ldots,c_n)=2^{-\lceil\frac{n}{2}\rceil}
\sum_{\varepsilon\in\{0,1\}^n}
\left(\prod_{1\leq j< n} (-1)^{\varepsilon_j\varepsilon_{j+1}}\right)
\left(\prod_{1\leq k\leq n} c_k^{\varepsilon_k}\right) \]
such that the order of the symbols $c_k$ is preserved in the products.
\end{point}
\begin{point}\label{po:classes}
Let $\overline {\mathbb C}=\mathbb C\cup\{\infty\}$ denote the Riemann sphere.
Some subsets are:
$\mathrm i\overline{\mathbb R}=
\mathrm i\mathbb R\cup\{\infty\}$, $\overline{\mathbb R^-}=(-\infty,0]\cup\{\infty\}$,
 $\overline{\mathbb C^-}=\{s\in\mathbb C\,:\,\Rea s\leq 0\}\cup\{\infty\}$,
${\mathring{\mathbb D}}^1=\{z\in\mathbb C\,:\,|z|< 1\}$.
We define the  functions
$\pol J=-\mathrm i \sgn \mathrm iJ$,
$|J|_{\mathrm i}=\sqrt{-J^2}$,
$|Q|_{\mathrm r}=\sqrt{Q^2}$,
$|P|_{\mathcal F}=\frac12-\left|\frac12-P\right|_{\mathrm r}$.
We have the following commutative diagram on certain subsets of the complex plane:
\begin{equation}
\xymatrix{(\arabic{equation})
&&J\in \mathbb C\setminus\overline{\mathbb R}\ar@/^/[r]^{J\mapsto\pol J}
\ar[dl]_{J\mapsto |J|_{\mathrm i}}\ar@/_1pc/[dll]_{J\mapsto -J^2}
&J\in\{\mathrm i,-\mathrm i\}\\
S\in\mathbb C\setminus\overline{\mathbb R^-}\ar@/^/[r]^{\,\,S\mapsto\sqrt S}\ar@/^/[d]^{S\mapsto \frac{1-S}4}
&Q\in \mathbb C\setminus\overline{\mathbb C^-}\ar[l]^{
\substack{\phantom{W}\\[-1mm]Q\mapsto Q^2}}\ar@/^/[d]^{Q\mapsto \frac{1-Q}2}&
Q\in\mathbb C\setminus \mathrm i\overline{\mathbb R}\ar@/^0.6pc/[ll]
\ar@/^/[r]^{Q\mapsto\sgn Q}\ar@/_/[l]_{\qquad Q\mapsto |Q|_{\mathrm r}}
\ar@/^/[d]^{Q\mapsto \frac{1-Q}2}
&Q\in\{1,-1\}\ar@/^/[d]^{Q\mapsto \frac{1-Q}2}\\
T\in\mathbb C\setminus\bigl(\frac14-\overline{\mathbb R^-}\bigr)
\ar@/^/[r]^{T\mapsto\sqrt[\mathcal F]T}\ar@/^/[u]^{T\mapsto 1-4T}&
P\in \mathbb C\setminus\bigl(\frac12-\overline{\mathbb C^-}\bigr)
\ar[l]^{\substack{\phantom{W}\\[0.5mm]\,\,P\mapsto P(1-P)}}\ar@/^/[d]^{P\mapsto -P(1-P)^{-1}}\ar@/^/[u]^{P\mapsto 1-2P}&
P\in\mathbb C\setminus\left(\frac12-\mathrm i\overline{\mathbb R}\right)\ar@/^0.7pc/[ll]
\ar@/^/[r]^{\qquad P\mapsto\idem P}\ar@/_/[l]_{P\mapsto |P|_{\mathcal F}}\ar@/^/[u]^{P\mapsto 1-2P}
&P\in\{0,1\}\ar@/^/[u]^{P\mapsto 1-2P}\\
&W\in {\mathring{\mathbb D}}^1\ar@/^/[u]^{
\substack{\phantom{W}\\[1mm]W\mapsto-W(1-W)^{-1}}}&&}\label{e5}\notag\addtocounter{equation}{1}\end{equation}
such that $\pol J$, $\sgn Q$, $\idem P$, $|J|_{\mathrm i}$, $|Q|_{\mathrm r}$, $|P|_{\mathcal F}$
yield idempotent operations and they yield decompositions
\[J=|J|_{\mathrm i}\pol J,\qquad
Q=|Q|_{\mathrm r}\sgn Q,\qquad
P=\idem P+|P|_{\mathcal F}-2|P|_{\mathcal F}\idem P.\]
\end{point}
\begin{defin}
Suppose that $\mathfrak A$ is a locally convex algebra and $R$ is a compact subset of $\overline{\mathbb C}$.
We define
$\spec_R(\mathfrak A)$
as the set containing all elements $X\in\mathfrak A$ such that the functions
\[f_R:\,z\in R\setminus\{\infty\}\mapsto (z-X)^{-1}\qand g_R:\,z\in R
\setminus\{0\}\mapsto X(1-z^{-1}X)^{-1}\]
are well-defined and continuous.
The topology of $\spec_R(\mathfrak A)$ is induced from the compact-open topology of the continuous functions
$(f_R)|_{\mathbb D^1\cap R}$ and $(g_R)|_{R\setminus{\mathring{\mathbb D}}^1}$.
(If $R$ is symmetric for conjugation then we may use the functions
$f_R(z)f_R(\bar z)$ and $g_R(z)g_R(\bar z)$ in order to get a formally real characterization.)
\end{defin}
Classes  of interest are like
$\spec_{\mathrm i\overline{\mathbb R}}(\mathfrak A)$, etc.,
i.~e.~the elements spectrally avoiding $\mathrm i\overline{\mathbb R}$, etc.
Another way to specify spectral conditions is to ask for
skew-involutions, involutions, or idempotents.
The main spectral classes correspond to the sets in \eqref{e5} for $\mathfrak A=\mathbb C$.
\begin{defin}
We define the corresponding formal spectral classes by
\begin{align}
J\in\widetilde\spec_{\overline{\mathbb R}}(\mathfrak A)&\Leftrightarrow
\text{ $\textstyle{\frac{1}{z}\left(\left(\frac{1+z}{2}\right)^2+
\left(\frac{1-z}{2}\right)^2J^2\right)}$
is invertible in $\mathfrak A[z,z^{-1}]$},
\notag\\
J\in\Skvol(\mathfrak A)&\Leftrightarrow
\text{ $J^2=-1$},
\notag\\
Q\in\widetilde\spec_{\mathrm i\overline{\mathbb R}}(\mathfrak A)&\Leftrightarrow
\text{ $\textstyle{\frac{1+z}2+\frac{1-z}2Q=\Lambda(z,Q)}$
is invertible in $\mathfrak A[z,z^{-1}]$},
\notag\\
Q\in\widetilde\spec_{\overline{\mathbb C^-}}(\mathfrak A)&\Leftrightarrow
\text{ $\textstyle{\frac{1+z}2+\frac{1-z}2Q=\Lambda(z,Q)}$ is invertible in $\mathfrak A[z]$},
\notag\\
Q\in\Invol(\mathfrak A)&\Leftrightarrow
\text{ $Q^2=1$},
\notag\\
S\in\widetilde\spec_{\overline{\mathbb R^-}}(\mathfrak A)&\Leftrightarrow
\text{ $\textstyle{\frac{1}{z}\left(\left(\frac{1+z}{2}\right)^2-
\left(\frac{1-z}{2}\right)^2S\right)=\Lambda(z,S,z^{-1})}$ is invertible in
$\mathfrak A[z,z^{-1}]$ }
\notag\\
P\in\widetilde\spec_{\frac12+\mathrm i\overline{\mathbb R}}(\mathfrak A)&\Leftrightarrow
\text{ $\textstyle{(1-P)+Pz}$
is invertible in $\mathfrak A[z,z^{-1}]$},
\notag\\
P\in\widetilde\spec_{\frac12-\overline{\mathbb C^-}}(\mathfrak A)&\Leftrightarrow
\text{ $\textstyle{(1-P)+Pz}$
is invertible in $\mathfrak A[z]$},
\notag\\
P\in\Idem(\mathfrak A)&\Leftrightarrow
\text{ $P^2=P$},
\notag\\
T\in\widetilde\spec_{\frac14-\overline{\mathbb R^-}}(\mathfrak A)&\Leftrightarrow
\text{ $\textstyle{1-(1-z)(1-z^{-1})T}$ is invertible in
$\mathfrak A[z,z^{-1}]$},
\notag\\
W\in\widetilde\spec_{\overline{\mathbb C}\setminus\mathring{\mathbb D}^1}(\mathfrak A)&\Leftrightarrow
\text{ $1+zW$ is invertible in $\mathfrak A[z]$}
.\notag
\end{align}
\end{defin}
\begin{point}
If $\mathfrak A$ is a locally convex algebra, then
the formal spectral classes and their ordinary counterparts are the same.
Indeed, the continuity of the resolvent terms implies smoothness by the resolvent identity,
hence the existence of the appropriate Fourier series,
and, conversely, the existence of the expansions implies continuity.
\end{point}
\begin{objective}\label{ob:our}
We want to establish the spectral correspondences and decompositions
as in point \ref{po:classes} for the formal spectral classes.
\end{objective}
\section{Calculations with $\frac12$}
\subsection*{\textbf{A. Sign and square root}}
\begin{defin} For $Q\in\widetilde\spec_{\mathrm i\overline{\mathbb R}}(\mathfrak A)$, we define
\[\sgn Q=\int \frac{\frac{1-z}{2}+\frac{1+z}{2}Q}{\frac{1+z}{2}+\frac{1-z}{2}Q}\,\frac{|dz|}{2\pi}=
\int\frac{\Lambda(-z,Q)}{\Lambda(z,Q)}\frac{|dz|}{2\pi}.\]
\end{defin}
\begin{prop}\label{lem:signprop}
If $Q\in\widetilde\spec_{\mathrm i\overline{\mathbb R}}(\mathfrak A)$,
then $-Q,Q^{-1}\in\widetilde\spec_{\mathrm i\overline{\mathbb R}}(\mathfrak A)$.
$Q$ commutes with $\sgn Q$. $\sgn -Q=-\sgn Q$ and $\sgn Q^{-1}=\sgn Q$. Moreover,
\[(\sgn Q)^2=1.\]
\begin{proof}
Substituting $z=-1$ we see that $Q^{-1}$ exists.
The first statement follows from the identities $\Lambda(z,-Q)=z\Lambda(z^{-1},Q)$
and $\Lambda(z,Q^{-1})=Q^{-1}\Lambda(-z,Q)$.
Furthermore, $Q$ and $\sgn Q$ commute, because $Q$ commutes with the integrand in $\sgn Q$.
The identities
\[\frac{\Lambda(-z,-Q)}{\Lambda(z,-Q)}\,\frac{|dz|}{2\pi}=
-\frac{\Lambda(-z^{-1},Q)}{\Lambda(z^{-1},Q)}\,\frac{|d(z^{-1})|}{2\pi}
\text{\quad and \quad}
\frac{\Lambda(-z,Q^{-1})}{\Lambda(z,Q^{-1})}\,\frac{|dz|}{2\pi}=
\frac{\Lambda(z,Q)}{\Lambda(-z,Q)}\,\frac{|d(-z)|}{2\pi}\]
integrated prove the first and second equalities, respectively.

The critical one is the involution property. We give several proofs.

\textit{``Matrix algebraic'' proof.}
Let $\mathsf H_{1/2}=\sum_{s\in \mathbb Z+\frac12}(\sgn s)\mathbf e_{s,s}$ be the
$\left(\mathbb Z+\frac12\right)\times \left(\mathbb Z+\frac12\right)$  matrix of the odd Hilbert transform.
Let $\mathsf U\left(\tfrac{1+Q}2z^{-1/2}+\tfrac{1-Q}2z^{1/2}\right)=
\sum_{s\in \mathbb Z+\frac12}\tfrac{1+Q}2\mathbf e_{s-\frac12,s}-\tfrac{1-Q}2\mathbf e_{s+\frac12,s}$ be the
$\mathbb Z\times\left(\mathbb Z+\frac12\right)$ matrix of the action of multiplication  by
$\left(\frac{1+Q}2+\frac{1-Q}2z\right)z^{-1/2}$. According to our assumption, this has a
$\left(\mathbb Z+\frac12\right)\times\mathbb Z$ inverse matrix representing the action of multiplication by
$z^{1/2}\left(\frac{1+Q}2+\frac{1-Q}2z\right)^{-1}$. So, we can consider the matrix
\[\mathsf B\left(\tfrac{1+Q}2z^{-1/2}+\tfrac{1-Q}2z^{1/2}\right)=
\mathsf U\left(\tfrac{1+Q}2z^{-1/2}+\tfrac{1-Q}2z^{1/2}\right)
\mathsf H_{1/2}
\mathsf U\left(\tfrac{1+Q}2z^{-1/2}+\tfrac{1-Q}2z^{1/2}\right)^{-1}.\]
Due to the special shape of the matrices involved, it is easy to see that this is an involution
which is the same as the even Hilbert transform $\mathsf H=\sum_{s\in\mathbb Z}(\sgn s)\mathbf e_{s,s}$,
except in the $0$th column.
This special shape implies that the diagonal element in the $0$th column
is an involution. On the other hand, it is easy to see that this diagonal element is exactly $\sgn Q$.

\textit{``Resolvent algebraic'' proof.}
As
\begin{equation}
1-(\sgn Q)^2=\iint\left(1-\frac{\Lambda(-z,Q)}{\Lambda(z,Q)}\frac{\Lambda(-w,Q)}{\Lambda(w,Q)}\right)\,
\frac{|dz|}{2\pi}\frac{|dw|}{2\pi},
\label{e6}\end{equation}
we should show that this integral is $0$.
This, however, follows from the key identity
\begin{multline}
1-\frac{\Lambda(-z,Q)}{\Lambda(z,Q)}\frac{\Lambda(-w,Q)}{\Lambda(w,Q)}=
\left(\frac{z+w}2\right)\frac{1-Q^2}{\Lambda(z,Q)\Lambda(w,Q)}
=\\=\frac12\left[\frac{z+w}{z-w}\right]\frac{(z-w)(1-Q^2)}{\Lambda(z,Q)\Lambda(w,Q)}
=\frac12\left[\frac{z+w}{z-w}\right]\left(\frac{(z-1)(1-Q^2)}{\Lambda(z,Q)}
-\frac{(w-1)(1-Q^2)}{\Lambda(w,Q)}\right),\label{e7}
\end{multline}
which does make sense in $\mathfrak A[[z^{-1},z]][[w^{-1},w]]$.
Indeed, \eqref{e6} can be continued as
\[=\iint  \left[\frac{z+w}{z-w}\right]\frac{(z-1)(1-Q^2)}{2\Lambda(z,Q)} \frac{|dz|}{2\pi}\frac{|dw|}{2\pi}
-\iint  \left[\frac{z+w}{z-w}\right]\frac{(w-1)(1-Q^2)}{2\Lambda(w,Q)} \frac{|dz|}{2\pi}\frac{|dw|}{2\pi}.\]
Evaluating the integrals we find
\[=\iint  0\cdot\frac{(z-1)(1-Q^2)}{2\Lambda(z,Q)} \frac{|dz|}{2\pi}-
\iint  0\cdot\frac{(w-1)(1-Q^2)}{2\Lambda(w,Q)} \frac{|dw|}{2\pi}=0-0=0.\]
This proof, like the previous one, relies heavily on the nature of Laurent series.
Nevertheless the argument can be modified so that formally it makes sense
in the analytical and the algebraic cases as well.

\textit{``Resolvent analytic'' proof.} According to the discussion about transformation kernels,
\eqref{e6} can be continued as follows:
\[=\lim_{t\nearrow1}\iint\tilde{\mathcal R}(t,wz^{-1})\left(1-
\frac{\Lambda(-z,Q)}{\Lambda(z,Q)}\frac{\Lambda(-w,Q)}{\Lambda(w,Q)}\right)\,
\frac{|dz|}{2\pi}\frac{|dw|}{2\pi}.\]
By simple arithmetic in the integrand, this yields
\[=\lim_{t\nearrow1}\biggl(\iint\tilde{\mathcal H}(t,wz^{-1})
\frac{(z-1)(1-Q^2)}{2\Lambda(z,Q)}\frac{|dz|}{2\pi}\frac{|dw|}{2\pi}
-\iint\tilde{\mathcal H}(t,wz^{-1})\frac{(w-1)(1-Q^2)}{2\Lambda(w,Q)}
\frac{|dz|}{2\pi}\frac{|dw|}{2\pi}\biggr).\notag\]
Executing the integrals we find
\[=\lim_{t\nearrow1}\biggl(\int0\cdot\frac{(z-1)(1-Q^2)}{\Lambda(z,Q)}\frac{|dz|}{2\pi}
-\int0\cdot\frac{(w-1)(1-Q^2)}{\Lambda(w,Q)}\frac{|dw|}{2\pi}\biggr)=0+0=0,\]
yielding, ultimately, the identity. We remark that in the analytic case, it would actually be simpler to
use the kernel $\mathcal R(t,wz^{-1})$.
\end{proof}
\end{prop}
\begin{defin}
If $S\in\widetilde\spec_{\overline{\mathbb R^-}}(\mathfrak A)$, then
we define the inverse square root operation as
\[\sqrt S=
\int \frac{zS}{\left(\frac{1+z}{2}\right)^2-\left(\frac{1-z}{2}\right)^2S}\,\frac{|dz|}{2\pi}
=\int\frac{S}{\Lambda(z,S,z^{-1})}\frac{|dz|}{2\pi}.\]
\end{defin}
\begin{prop}
Suppose that $S\in\widetilde\spec_{\overline{\mathbb R^-}}(\mathfrak A)$.
Then $S^{-1}\in\widetilde\spec_{\overline{\mathbb R^-}}(\mathfrak A)$.
The elements $\sqrt S$ and $S$ commute with each other.
$\sqrt{S^{-1}}={\sqrt S}^{-1}$. Furthermore,
\begin{equation}(\sqrt S)^2=S.\label{e:sq1} \end{equation}
\begin{proof}
Substituting $z=-1$ into the resolvent term, we see that $S^{-1}$ exists. The identity
$\dfrac1{\Lambda(z,S^{-1},z^{-1})}=\dfrac S{\Lambda(-z,S,(-z)^{-1})}$
shows that $S^{-1}\in\widetilde\spec_{\overline{\mathbb R^-}}(\mathfrak A)$.
Integrated, it yields $S\sqrt{S^{-1}}=\sqrt S$.
If the square-root identity  \eqref{e:sq1} holds, then this implies $\sqrt{S^{-1}}=\sqrt S^{-1}$.
So, it remains to prove \eqref{e:sq1}.
As
\begin{equation}(\sqrt S)^2-S=\iint\left(\frac{S}{\Lambda(z,S,z^{-1})}
\frac{S}{\Lambda(w,S,w^{-1})} -S\right)\,\frac{|dz|}{2\pi}\frac{|dw|}{2\pi},
\label{e8}\end{equation}
we have to show that this integral is $0$.
This follows using the key identities
\begin{multline}
\frac{S}{\Lambda(z,S,z^{-1})}\frac{S}{\Lambda(w,S,w^{-1})}-S=\\=
S\frac{S-\Lambda(z,S,z^{-1},1,w,S,w^{-1})}{2\Lambda(z,S,z^{-1})\Lambda(w,S,w^{-1})}
+S\frac{S-\Lambda(z,S,z^{-1},1,w^{-1},S,w)}{2\Lambda(z,S,z^{-1})\Lambda(w,S,w^{-1})};
\label{e9}\end{multline}
\begin{equation}\frac{S-\Lambda(z,S,z^{-1},1,w,S,w^{-1})}{2\Lambda(z,S,z^{-1})\Lambda(w,S,w^{-1})}=
\frac12\left[\frac{zw^{-1}+1}{zw^{-1}-1}\right]
\left(\frac{\Lambda(-z,S,z^{-1})}{\Lambda(z,S,z^{-1})}-\frac{\Lambda(-w,S,w^{-1})}{\Lambda(w,S,w^{-1})}
\right);\label{e10}\end{equation}
\begin{equation}\frac{S-\Lambda(z,S,z^{-1},1,w^{-1},S,w)}{2\Lambda(z,S,z^{-1})\Lambda(w,S,w^{-1})}=
\frac12\left[\frac{zw+1}{zw-1}\right]
\left(\frac{\Lambda(-z,S,z^{-1})}{\Lambda(z,S,z^{-1})}-\frac{\Lambda(w,S,-w^{-1})}{\Lambda(w,S,w^{-1})}
\right)\label{e11}.\end{equation}
Indeed, after we decomposed the integrand in \eqref{e8} according to (\ref{e9}--\ref{e11}),
we can show that both parts are $0$
as we did in the previous proof.
\end{proof}
\end{prop}
\begin{prop}\label{lem:rootsign}
$Q\in\widetilde\spec_{\mathrm i\overline{\mathbb R}}(\mathfrak A)$ if and only if
$Q^2\in\widetilde\spec_{\overline{\mathbb R^-}}(\mathfrak A)$. In this case
\[\sgn Q=Q^{-1}\sqrt{Q^2}.\]
\begin{proof}
The first statement follows from the equality $\Lambda(z,Q^2,z^{-1})=\Lambda(z,Q)\Lambda(z^{-1},Q)$.
The identity statement follows from
\begin{multline}\sgn Q=\frac{\sgn Q}{2}+\frac{\sgn Q}{2}=
\frac{1}{2}\int\frac{\Lambda(-z,Q)}{\Lambda(z,Q)}\,\frac{|dz|}{2\pi}+
\frac{1}{2}\int \frac{\Lambda(-z^{-1},Q)}{\Lambda(z^{-1},Q)}
\,\frac{|d(z^{-1})|}{2\pi} \\
=\int\frac{1}{2}\left(\frac{\Lambda(-z,Q)}{\Lambda(z,Q)}+ \frac{\Lambda(-z^{-1},Q)}{\Lambda(z^{-1},Q)}
\right)\,\frac{|dz|}{2\pi}=\int\frac{Q}{\Lambda(z,Q^2,z^{-1})}
\,\frac{|dz|}{2\pi}=Q^{-1}\sqrt{Q^2}.\qquad\notag\end{multline}
\par \vspace{-\baselineskip}\end{proof}
\end{prop}

\subsection*{\textbf{B. Finer analysis of the resolvent terms}}
\begin{defin} (a) We define $|Q|_{\mathrm r}=\sqrt{Q^2}=Q\sgn Q$.

(b) If $F\in\mathfrak A$ is an involution, then an element $A\in\mathfrak A$ can be written in matrix form
\[\begin{bmatrix}\frac{1-F}2A\frac{1-F}2&\frac{1-F}2A\frac{1+F}2\\
\frac{1+F}2A\frac{1-F}2&\frac{1+F}2A\frac{1+F}2\end{bmatrix}.\]

Suppose that $Q\in\widetilde\spec_{\mathrm i\overline{\mathbb R}}(\mathfrak A)$.
In the decomposition of $\mathfrak A$ along the involution $\sgn Q$, the various components are denoted according to
\[\sgn Q=\begin{bmatrix} -1_{-\sgn Q} &\\& 1_{\sgn Q}\end{bmatrix},\qquad Q=
\begin{bmatrix} -Q^- &\\& Q^+\end{bmatrix},\qquad |Q|_{\mathrm r}=\begin{bmatrix} Q^-&\\& Q^+\end{bmatrix}.\]
\end{defin}
\begin{prop} $Q^\pm\in\widetilde\spec_{\mathrm i\overline{\mathbb R}}(\mathfrak A1_{\pm\sgn Q})$,
$\sgn Q^\pm=1_{\pm\sgn Q}$, and $\sgn |Q|_{\mathrm r}=1$.
\begin{proof}
The decomposition of $Q$ along $\sgn Q$ allows us to consider $Q$
separately in the direct sum components of $\mathfrak A$.
In particular, the sign integral splits, too, and it necessarily yields
$\sgn Q^+=1_{\sgn Q}$ and $\sgn -Q^-=-1_{\sgn -Q}$.
 The statement follows from this immediately.
\end{proof}
\end{prop}
\begin{prop}\label{lem:decoinv}
If $Q\in\widetilde\spec_{\mathrm i\overline{\mathbb R}}(\mathfrak A)$, then

(a)
\[\frac{1}{\frac{1+z}{2}+\frac{1-z}{2}Q}=\frac{1}{\Lambda(z,Q)}\]
is given by
\begin{multline}
\ldots+
\begin{bmatrix}\left(\frac{Q^--1}{Q^-+1}\right)^2\frac{2}{Q^-+1}&\\&0\end{bmatrix}z^{-3}+
\begin{bmatrix}\left(\frac{Q^--1}{Q^-+1}\right)\frac{2}{Q^-+1}&\\&0\end{bmatrix}z^{-2}+
\begin{bmatrix}\frac{2}{Q^-+1}&\\&0\end{bmatrix}z^{-1}+
\\+\begin{bmatrix} 0&\\&\frac{2}{Q^++1}  \end{bmatrix}1+
\begin{bmatrix}0&\\&\frac{2}{Q^++1}\left(\frac{Q^+-1}{Q^++1}\right)\end{bmatrix}z+
\begin{bmatrix}0&\\&\frac{2}{Q^++1}\left(\frac{Q^+-1}{Q^++1}\right)^2 \end{bmatrix}z^2+\ldots\notag
\end{multline}

(b)
\[\frac{\frac{1-z}{2}+\frac{1+z}{2}Q}{\frac{1+z}{2}+\frac{1-z}{2}Q}=\frac{\Lambda(-z,Q)}{\Lambda(z,Q)}\]
is given by
\[
\ldots+\begin{bmatrix}  -2\left(\frac{Q^--1}{Q^-+1}\right)^3&\\&0\end{bmatrix}z^{-3}+
\begin{bmatrix}  -2\left(\frac{Q^--1}{Q^-+1}\right)^2&\\&0\end{bmatrix}z^{-2}+
\begin{bmatrix}  -2\left(\frac{Q^--1}{Q^-+1}\right)&\\&0\end{bmatrix}z^{-1}+
\]\[
+\begin{bmatrix} -1_{-\sgn Q}&\\& 1_{\sgn Q}\end{bmatrix}1
+\begin{bmatrix}0&\\&2\left(\frac{Q^+-1}{Q^++1}\right) \end{bmatrix}z+
\begin{bmatrix}0&\\&2\left(\frac{Q^+-1}{Q^++1}\right)^2  \end{bmatrix}z^2+
\begin{bmatrix}0&\\&2\left(\frac{Q^+-1}{Q^++1}\right)^3 \end{bmatrix}z^3+\ldots\notag
\]

(c) In particular,
\[\frac1{1+|Q|_{\mathrm r}}=\int\frac12\frac{1+z}{\Lambda(z,Q)}\frac{|dz|}{2\pi}.\]
\begin{proof}
(a) It is enough to consider the $Q^+$ part of the decomposition, because the other part follows
from changing $z$ to $z^{-1}$.  So, we can suppose that $Q=|Q|$ and $\sgn Q=1$.

We try to figure out the coefficients in the expansion $\frac{1}{\Lambda(z,Q)}=\sum_{n\in\mathbb Z}a_nz^n.$
The equality $\sgn Q=1$ means that
\begin{equation}
0=\frac{1-\sgn Q}{2}=\frac{1}{2}\int\left(1-\frac{\Lambda(-z,Q)}{\Lambda(z,Q)}\right)\frac{|dz|}{2\pi}
=\int\frac{1-Q}{2}\frac{z}{\Lambda(z,Q)}\frac{|dz|}{2\pi}=\frac{1-Q}{2}a_{-1}.\label{e12}\end{equation}
The  product of $\Lambda(z,Q)=\frac{1+Q}{2}+\frac{1-Q}{2}z$ and $\sum_{n\in\mathbb Z}a_nz^n$ gives $1$, so
\begin{equation}\frac{1-Q}{2}a_{n-1}+\frac{1+Q}{2}a_{n}=\boldsymbol\delta_{0,n}1\label{e13}.\end{equation}
From \eqref{e12} and the case $n=0$ in \eqref{e13}, we obtain that $a_0=\frac{2}{1+Q}$.
After that, from \eqref{e13}, we find $a_{n+1}=\tfrac{Q-1}{Q+1}a_n$, yielding the positive-numbered coefficients.
Then $F(z)=\sum_{n\in\mathbb N}a_nz^n$
already inverts $\Lambda(z,Q)$, hence, from the uniqueness of the inverse, $a(z)=F(z)$.

(b) and (c) follow from part (a) by  simple algebra.
\end{proof}
\end{prop}
\begin{prop} \label{lem:sqex}
(a) If $S\in\widetilde\spec_{\overline{\mathbb R^-}}(\mathfrak A)$, then
\[\frac{z}{\left(\frac{1+z}{2}\right)^2-\left(\frac{1-z}{2}\right)^2S}=\frac{1}{\Lambda(z,S,z^{-1})}\]
yields the expansion
\[\sqrt S^{-1}1
+\sqrt S^{-1}\left(\frac{\sqrt S-1}{\sqrt S+1}\right)(z+z^{-1})+
\sqrt S^{-1}\left(\frac{\sqrt S-1}{\sqrt S+1}\right)^2(z^2+z^{-2})+\ldots\]

(b) In particular,
\[\frac1{\sqrt S+1}=\int\frac12\frac{1+z}{\Lambda(z,S,z^{-1})}\frac{|dz|}{2\pi}. \]
\begin{proof}
(a) Take $Q=\sqrt S$. Proposition \ref{lem:rootsign} yields $\sgn Q=1$.
Applying  the identity
$\dfrac{1}{\Lambda(z,S,z^{-1})}=\dfrac{1}{\Lambda(z,Q)\Lambda(z^{-1},Q)}$
and Proposition \ref{lem:decoinv}, it follows that the coefficient
of $z^{\pm n}$ ($n\geq0$) in the  expansion  is
\[\left(\frac{2}{Q+1}\right)^2\sum_{m\in\mathbb N}\left(\frac{Q-1}{Q+1}\right)^{n+2m}=
\left(\frac{2}{Q+1}\right)^2\left(\frac{Q-1}{Q+1}\right)^{n}
\left(1-\left(\frac{Q-1}{Q+1}\right)^2\right)^{-1},\]
which simplifies as above.
The rapid decrease of $\left(\frac{Q-1}{Q+1}\right)^s$ makes our computations legal.

(b) follows from part (a).
\end{proof}
\end{prop}
\begin{prop}
(a) $Q\in\widetilde\spec_{\overline{\mathbb C^-}}(\mathfrak A)$
if and only if $Q\in\widetilde\spec_{\mathrm i\overline{\mathbb R^-}}(\mathfrak A)$ and $\sgn Q=1$.

(b) $W\in\widetilde\spec_{\overline{\mathbb C}\setminus\mathring{\mathbb D}^1}(\mathfrak A)$
if and only if $W^n$ is rapidly decreasing.

(c) If $Q\in\widetilde\spec_{\overline{\mathbb C^-}}(\mathfrak A)$,
then $\frac{1-Q}{1+Q}\in\widetilde\spec_{\overline{\mathbb C}\setminus\mathring{\mathbb D}^1}(\mathfrak A)$.
Conversely, if $W\in\widetilde\spec_{\overline{\mathbb C}\setminus\mathring{\mathbb D}^1}(\mathfrak A)$,
then $\frac{1-W}{1+W}\in\widetilde\spec_{\overline{\mathbb C^-}}(\mathfrak A)$.
This establishes a bijection.
\begin{proof}
(a) follows from Proposition \ref{lem:decoinv}.a.
(b) holds because in those cases $(1-Wz)^{-1}=\sum_{n\in\mathbb N}W^nz^n$ must hold.
(c) follows from
$\Lambda\left(z,\frac{1-W}{1+W}\right)=\frac{1+zW}{1+W}$ and
$1+\frac{1-Q}{1+Q}z=\frac{2\Lambda(z,Q)}{1+Q}$.
\end{proof}
\end{prop}
\subsection*{\textbf{C. On our objective}}
\begin{point}
Now, it is easy to see that the propositions proven above are sufficient to establish
all the spectral correspondences asked in \ref{ob:our}.
We merely define $|J|_{\mathrm i}=\sqrt{-J^2}$, $\pol J=J \,|-J^2|_{\mathrm i}^{-1}$,
$\idem P=\frac12-\frac12\sgn\left(\frac12-P\right)$,
$|P|_{\mathcal F}=\frac12-\left|\frac12-P\right|_{\mathrm r}$,
and $\sqrt[\mathcal F] T=\frac12-\sqrt{\frac14-T}$.
Hence our objective is established.
\end{point}
This is, however, not to say that everything is just like for locally convex algebras:

\subsection*{\textbf{D. Comparison to the case of locally convex algebras}}

\begin{prop}\label{lem:rootreal} If $\mathfrak A$ is a locally convex algebra, then
the condition that $S\in\spec_{\overline{\mathbb R^-}}(\mathfrak A)$  is equivalent to the
condition that the function
\[\frac{1+t}{2}+\frac{1-t}{2}S\]
has a continuous inverse on $[-1,1]$. (This is the same thing as to say that the segment connecting
$1$ and $S$ is continuously invertible.)
The square root can  be expressed as
\[\sqrt S=\int_{t\in[-1,1]}\frac{S}{\frac{1+t}{2}+\frac{1-t}{2}S}\,\frac{dt}{\pi\sqrt{1-t^2}}.\]
\begin{proof}
It follows by change of variables using
$t=\frac{z+z^{-1}}{2}$.
\end{proof}
\end{prop}
\begin{point}
If $S\in\spec_{\overline{\mathbb C^-}}(\mathfrak A)$,
then $\frac{1+t}{2}+\frac{1-t}{2}S$ $(t\in[-1,1])$ is clearly invertible.
As a rapidly decreasing power series in $t$, using $t=\frac{z+z^{-1}}2$
and considering the coefficients of $z^k$ in  $\left(\frac{z+z^{-1}}2\right)^n$,
it follows that $S\in\spec_{\overline{\mathbb R^-}}(\mathfrak A)$.
Hence the inclusion
$\spec_{\overline{\mathbb C^-}}(\mathfrak A)\subset\spec_{\overline{\mathbb R^-}}(\mathfrak A)$
is true.
In the general context, this cannot be done so, because the boundedness of the elements
$\frac1{2^n}\left(\begin{smallmatrix}n \\ k\end{smallmatrix}\right)$ is not always clear.
Similar comment applies for
$\spec_{\frac12+\overline{\mathbb C^-}}(\mathfrak A)\subset\spec_{\frac12+\overline{\mathbb R^-}}(\mathfrak A)$.
\end{point}

\section{ Formal homotopies}
\begin{point} One expects certain natural behaviour from the operations above.
For example, one expects to have a homotopy from $Q$ to $\sgn Q$ inside
$\widetilde\spec_{\mathrm i\overline{\mathbb R}}(\mathfrak A)$.
In general algebras, one cannot use continuous variables, but one can
come up with  homotopies using formal variables.
Let us remind that an element
$Q\in\widetilde\spec_{\mathrm i\overline{\mathbb R}}(\mathfrak A)$
can be decomposed to a commuting pair, $\sgn Q$ and a perturbation of $1$ which is $|Q|_{\mathrm r}$.
But we may also consider this as a decomposition  to the commuting pair $\sgn Q$ and a perturbation of $0$ which is
\[\pert_{\mathrm r}Q=\frac{|Q|_{\mathrm r}-1}{|Q|_{\mathrm r}+1}
=-\int \frac{\frac{1+z}{2}-\frac{1-z}{2}Q}{\frac{1+z}{2}+\frac{1-z}{2}Q}\,\frac{|dz|}{2\pi}=
-\int\frac{\Lambda(z,-Q)}{\Lambda(z,Q)}\frac{|dz|}{2\pi}.\]
If we replace $\pert_{\mathrm r}Q$ by $t\pert_{\mathrm r}Q$ in the decomposition, then we obtain a
homotopy, this appears as $K(t,-1,Q)$ in what follows.
\end{point}
\begin{defin}
We define
\[K(t,z,Q)=\frac{1+\sgn Q}{2}\frac{\Lambda(tz,|Q|_{\mathrm r})}{\Lambda(t,|Q|_{\mathrm r})}+
\frac{1-\sgn Q}{2}z\frac{\Lambda(tz^{-1},|Q|_{\mathrm r})}{\Lambda(t,|Q|_{\mathrm r})},\]
\[H(t,z,Q)=\frac{1+\sgn Q}{2}\frac{\Lambda(t,|Q|_{\mathrm r})}{\Lambda(tz,|Q|_{\mathrm r})}+
\frac{1-\sgn Q}{2}z^{-1}\frac{\Lambda(t,|Q|_{\mathrm r})}{\Lambda(tz^{-1},|Q|_{\mathrm r})},\]
\[L(t,z,Q)=\frac{1+\sgn Q}{2}\Lambda(tz,|Q|_{\mathrm r})+\frac{1-\sgn Q}{2}z\Lambda(tz^{-1},|Q|_{\mathrm r}),\]
\[G(t,z,Q)=\frac{1+\sgn Q}{2}\frac{1}{\Lambda(tz,|Q|_{\mathrm r})}+
\frac{1-\sgn Q}{2}z^{-1}\frac{1}{\Lambda(tz^{-1},|Q|_{\mathrm r})}.\]
\end{defin}
\begin{prop}
The expressions $K(t,z,Q)$ and $H(t,z,Q)$ are multiplicative inverses of each other.
\[K(t,z,Q)=\frac{\Lambda(z,\sgn Q,t,|Q|_{\mathrm r})}{\Lambda(t,|Q|_{\mathrm r})}
=\Lambda\left(z,\frac{\Lambda(-t,|Q|_{\mathrm r})}{\Lambda(t,|Q|_{\mathrm r})}\sgn Q\right).\]
$K(t,1,Q)=1$, $K(1,-1,Q)=Q$, $K(0,-1,Q)=\sgn Q$, $K(-1,-1,Q)=Q^{-1}$.
\[H(t,z,Q)=\frac{\Lambda(t,|Q|_{\mathrm r})
\Lambda(z^{-1},\sgn Q,t,|Q|_{\mathrm r})}{\Lambda(tz,|Q|_{\mathrm r})\Lambda(tz^{-1},|Q|_{\mathrm r})}.\]
Similarly, the expressions $L(t,z,Q)$ and $G(t,z,Q)$ are inverses.
\[L(t,z,Q)=\Lambda(z,\sgn Q,t,|Q|_{\mathrm r})=\Lambda(t,\sgn Q,z,Q).\]
$L(t,1,Q)=\Lambda(t,|Q|_{\mathrm r})$, $L(1,-1,Q)=Q$, $L(0,-1,Q)=\frac12(Q+\sgn Q)$, $L(-1,-1,Q)=\sgn Q $.
\[G(t,z,Q)=\frac{\Lambda(z^{-1},\sgn Q,t,|Q|_{\mathrm r})}{
\Lambda(tz,|Q|_{\mathrm r})\Lambda(tz^{-1},|Q|_{\mathrm r})}.\]
\begin{proof}
The computation is easy if we notice that in the defining formulas the
coefficients of $\frac{1+\sgn Q}{2}$ and $\frac{1-\sgn Q}{2}$ live separate
lives because $\sgn Q$ is an involution.
\end{proof}
\end{prop}
The properties of $K$ show that $\frac{\Lambda(-t,|Q|_{\mathrm r})}{\Lambda(t,|Q|_{\mathrm r})}\sgn Q$ is
a homotopy from
$Q$ $(t=1)$ to $\sgn Q$ $(t=0)$ inside $\widetilde\spec_{\mathrm i\overline{\mathbb R}}(\mathfrak A)$.
Here the meaning of ``inside'' is  that the whole expression
satisfies the appropriate formal spectral condition.
\begin{prop}
\[K(t,w,Q)=\int \mathcal P(t,z)\frac{\Lambda(zw,Q)}{\Lambda(z,Q)}\,\frac{|dz|}{2\pi},\]
\[H(t,w,Q)=\int \mathcal P(t,z)\frac{\Lambda(z,Q)}{\Lambda(zw,Q)}\,\frac{|dz|}{2\pi},\]
\[G(t,w,Q)=\int \mathcal L(t,z)\frac1{\Lambda(zw,Q)}\frac{|dz|}{2\pi}.\]
\begin{proof}
This follows from the series expansion in  Proposition \ref{lem:decoinv}.
\end{proof}
\end{prop}
\begin{remark} In locally convex algebras,
 the controllability of the powers of $\frac{z+z^{-1}}2$
makes possible to consider
\[\Lambda\left(z,\Lambda(t,|Q|_{\mathrm r})\sgn Q\right)=
\frac{1+\sgn Q}2\Lambda(\Lambda(z,t),|Q|_{\mathrm r})+\frac{1-\sgn Q}2z\Lambda(\Lambda(z^{-1},t),|Q|_{\mathrm r}),\]
whose inverse turns out to be
\[\frac{\Lambda(z^{-1},\Lambda(t,|Q|_{\mathrm r})\sgn Q)
}{\Lambda(\Lambda(z,t),|Q|_{\mathrm r})\Lambda(\Lambda(z^{-1},t),|Q|_{\mathrm r})}.\]
This shows that  $\Lambda(t,|Q|_{\mathrm r})\sgn Q=\frac{1+t}2\sgn Q+\frac{1-t}2Q$ is
also a formal homotopy between $Q$ $(t=-1)$ and $\sgn Q$ $(t=1)$ inside
$\widetilde\spec_{\mathrm i\overline{\mathbb R}}(\mathfrak A)$.
\end{remark}
\begin{point} Similarly, we can contract elements inside
$\widetilde\spec_{\overline{\mathbb R^-}}(\mathfrak A)$  to $1$.
For $S\in\widetilde\spec_{\overline{\mathbb R^-}}(\mathfrak A)$ consider
\[C(t,S)=\left(\frac{1+t\frac{\sqrt S-1}{\sqrt S+1}}{1-t\frac{\sqrt S-1}{\sqrt S+1}}\right)^2,
\qquad
\sqrt{C(t,S)}=\frac{1+t\frac{\sqrt S-1}{\sqrt S+1}}{1-t\frac{\sqrt S-1}{\sqrt S+1}}.\]
The substitution $t\mapsto -t$ inverts them multiplicatively. In fact, the corresponding loops invert:
\end{point}
\begin{prop}For $S\in\widetilde\spec_{\overline{\mathbb R^-}}(\mathfrak A)$, we have
\[\frac1{\Lambda(w,C(t,S),w^{-1})}=\frac{\sqrt S}{\sqrt{C(t,S)}}\int \mathcal P(t,z)
 \frac1{\Lambda(zw,S,(zw)^{-1})}\frac{|dz|}{2\pi}. \]
\end{prop}
\begin{remark} In locally convex algebras, alternative contracting paths are rather trivial to
find. It is more interesting to see that the class of loops of type
\[ \frac1{\frac{1}{z}\left(\left(\frac{1+z}{2}\right)^2A^{-1}-
\left(\frac{1-z}{2}\right)^2B^{-1}\right)} \]
remains invariant with respect to the Poisson kernel. For $t=0$, they  contract to the geometric mean
\[\sqrt{A\cdot B}=\int  \frac1{\frac{1}{z}\left(\left(\frac{1+z}{2}\right)^2A^{-1}-
\left(\frac{1-z}{2}\right)^2B^{-1}\right)}\frac{|dz|}{2\pi}.\]
\end{remark}

\section{Calculations without $\frac12$}
As we have seen, much can be generalized to the case $\frac12\in\mathfrak A$.
It is natural to ask what happens in the lack of $\frac12$.
Then only the lower portion of \eqref{e5} can be generalized.
Again, the idempotent and the $\mathcal F$-square-root identities are the key properties.
\begin{defin} For $P\in\spec_{\frac12+\mathrm i\overline{\mathbb R}}(\mathfrak A)$, we define
\[\idem P=\int\frac{Pz}{1-P+Pz}\,\frac{|dz|}{2\pi}.\]
\end{defin}
\begin{prop}If $P\in\spec_{\frac12+\mathrm i\overline{\mathbb R}}(\mathfrak A)$,
then $1-2P$ is invertible,
$1-P,\frac{-P}{1-2P}\in\spec_{\frac12+\mathrm i\overline{\mathbb R}}(\mathfrak A)$, and
$\idem \, (1-P)=1-\idem P$, $\idem\frac{-P}{1-2P}=\idem P$. Furthermore,
\[(\idem P)^2=\idem P.\]
\begin{proof}
The invertibility statement follows from the substitution $z=-1$.
The identities $P+(1-P)z=(P+(1-P)z^{-1})z$ and
$(1-\frac{-P}{1-2P})+\frac{-P}{1-2P}z=(1-2P)^{-1}(1-P+P(-z))$
imply the spectral statements.
The identities
\[\frac{Pz}{(1-P)+Pz}\frac{|dz|}{2\pi}=\left(1-\frac{(1-P)z^{-1}}{P+(1-P)z^{-1}} \right)\frac{|d(z^{-1})|}{2\pi},\]
\[\frac{\frac{-P}{1-2P}z}{(1-\frac{-P}{1-2P})+\frac{-P}{1-2P}z}\frac{|dz|}{2\pi}=
\frac{P(-z)}{1-P+P(-z)} \frac{|d(-z)|}{2\pi}\]
integrated prove the first and second equalities, respectively.
We can prove the idempotent identity in several ways:

\textit{Matrix algebraic proof.} We can proceed as before, but have to conjugate
not the Hilbert transform involution but the idempotent  $\sum_{s\in-\mathbb N-\frac12}\mathbf e_{s,s}$.

\textit{A direct algebraic proof.} See Karoubi \cite{K}, Lemma III.1.23--24.

\textit{A resolvent algebraic proof.}
We should prove that
\begin{equation}\idem P(1-\idem P)=\iint\frac{Pz}{1-P+Pz}\frac{1-P}{1-P+Pw}\frac{|dz|}{2\pi}\frac{|dw|}{2\pi}
\label{e14}\end{equation}
is equal to $0$. It is natural try the proof along the steps
\begin{multline}
\frac{Pz}{1-P+Pz}\frac{1-P}{1-P+Pw}\sim \left[\frac{z+w}2\right]\frac{P(1-P)}{(1-P+Pz)(1-P+Pw)}=
\\=\left[\frac12\,\frac{z+w}{z-w}\right]\left(\frac{(z-1)P(1-P)}{1-P+Pz}-\frac{(w-1)P(1-P)}{1-P+Pw}\right)\sim 0,
\label{e15}\end{multline}
except it seems to be plagued by $\frac12$'s as before.
We have to demonstrate that the use of division by $2$ is of superficial nature in the proof.
This can be done as follows.

For a Laurent series $a(z,w)$, we define $:a(z,w):_{z,w}$ by linear extension from
\[:z^nw^m:_{z,w}=z^{\max(n,m)}w^{\min(n,m)}.\]
\begin{lemma} For any Laurent series $a(z,w)$, we have
\[\iint a(z,w)\frac{|dz|}{2\pi}\frac{|dw|}{2\pi}=\iint:a(z,w):_{z,w}\frac{|dz|}{2\pi}\frac{|dw|}{2\pi}.\]
\end{lemma}
For a symmetric Laurent series $a(z,w)$, i. e. such that $a(z,w)=a(w,z)$, we define formally
\[:\left[\frac{z+w}2\right]a(z,w):_{z,w}=:za(z,w):_{z,w}.\]

Suppose that $a(z,w)$ is anti-symmetric in its variables, included that the coefficient of $z^nw^n$ is always $0$.
We define $:\left[\frac12\,\frac{z+w}{z-w}\right]a(z,w):_{z,w}$
formally to be as it should be according to our natural expectations.
We have to check that the resulting expression is integral in terms of the coefficients.
In the present case, the definition yields by linear extension from
\[:\left[\frac12\,\frac{z+w}{z-w}\right](z^nw^m-z^mw^n):_{z,w}=z^nw^m+2\sum_{0<k<\frac{n-m}2} z^{n-k}w^{m+k}+
\boldsymbol\delta_{\frac{n+m}2\in\mathbb Z} z^{\frac{n+m}2}w^{\frac{n+m}2},\]
where $n>m$. Checking for elements of suitable bases it  is easy to see the following lemmas:
\begin{lemma} For any symmetric Laurent series $a(z,w)$,
\[:\left[\frac{z+w}2\right]a(z,w):_{z,w}= :\left[\frac12\,\frac{z+w}{z-w}\right] (z-w)a(z,w):_{z,w}.\]
\end{lemma}
\begin{lemma} For any Laurent series $a(z)$, we have
\[\iint :\left[\frac12\,\frac{z+w}{z-w}\right] (a(z)-a(w)):_{z,w}\frac{|dz|}{2\pi}\frac{|dw|}{2\pi}=0.\]
\end{lemma}
Now it is easy to carry out the proof.
From \eqref{e14} we should pass to the ``normal ordered'' form, after which the subsequent
manipulations as in \eqref{e15},  leading to $0$, make sense.
\end{proof}
\end{prop}
\begin{defin}
For $T\in\spec_{\frac14-\overline{\mathbb R^-}}(\mathfrak A)$, we define
\[\sqrt[\mathcal F] T=\int\frac{(1+z)T}{1+(z-2+z^{-1})T}\,\frac{|dz|}{2\pi}.\]
\end{defin}
\begin{defin}\label{def:div}
 Let $\mathfrak A\langle z\rangle^+$ be the space of formal Laurent series
$a(z)=a_0+\sum_{k=1}^\infty a_k\frac{z^k+z^{-k}}2$,
and let $\mathfrak A[z]^+$ be the space of formal Laurent series
$b(z)=b_0+\sum_{k=1}^\infty b_k(z^k+z^{-k})$.
Similarly,
let $\mathfrak A\langle z\rangle^-$ be the space of formal Laurent series
$a(z)=\sum_{k=1}^\infty a_{-k}\frac{z^k-z^{-k}}2$,
and let $\mathfrak A[z]^-$ be the space of formal Laurent series
$b(z)=\sum_{k=1}^\infty b_{-k}(z^k-z^{-k})$.
It is easy to see that
$\mathfrak A\langle z\rangle=\mathfrak A\langle z\rangle^+\oplus\mathfrak A\langle z\rangle^-$ is a natural
$\mathfrak A[z]^\pm=\mathfrak A[z]^+\oplus\mathfrak A[z]^-$-module,
in fact, this action is $\boldsymbol Z_2$-graded.
Multiplication of $1\in \mathfrak A\langle z\rangle$
yields a natural map from $\mathfrak A[z]^\pm$ into $\mathfrak A\langle z\rangle$:
\[1\cdot\left(\sum_{k=1}^\infty b_{-k}(z^k-z^{-k})+b_0+\sum_{k=1}^\infty b_k(z^k+z^{-k})\right)=
\sum_{k=1}^\infty 2b_{-k}\tfrac{z^k-z^{-k}}2+b_0+\sum_{k=1}^\infty 2b_k\tfrac{z^k+z^{-k}}2.\]
In fact, this notation can be extended to $b(z)\in\mathfrak A[z,z^{-1}]$ in a compatible way,
by $1\cdot z=\frac{z+z^{-1}}2+\frac{z-z^{-1}}2$, etc.
Integration can be defined for elements of $\mathfrak A\langle z\rangle$ or $\mathfrak A[z]^\pm$. Again,
it singles out the $0$th coefficient.
We see that if $a(z)\in\mathfrak A\langle z\rangle$ and
$b(z)\in\mathfrak A[z]^\pm$ as above, then
\[\int a(z)\cdot b(z)\frac{|dz|}{2\pi}=\sum_{n\in\mathbb Z} a_kb_k.\]
For example, $\int \left(1+\frac{z+z^{-1}}2\right)\cdot b(z)\frac{|dz|}{2\pi}=b_0+b_1$.
Furthermore, $\int b(z)\frac{|dz|}{2\pi}=\int 1\cdot b(z)\frac{|dz|}{2\pi}$.

We can extend this formalism to multiple variables.
The spaces $\mathfrak A\langle z\rangle\langle w\rangle$ and $\mathfrak A[z]^\pm[w]^\pm$
can be considered.
In the case of $\mathfrak A\langle z\rangle\langle w\rangle$,  colloquial notation  like
\[\tfrac{zw^{-1}+wz^{-1}}2\equiv\tfrac{z+z^{-1}}2\tfrac{w+w^{-1}}2-\tfrac{z-z^{-1}}2\tfrac{w-w^{-1}}2\]
is allowed. However, an other space between
$\mathfrak A\langle z\rangle\langle w\rangle$ and $\mathfrak A[z]^\pm[w]^\pm$
can be considered. Indeed, let  $\mathfrak A\{z,w\}$
be the space of the formal combinations of the basis elements
\[1, \tfrac{z^n+z^{-n}}2,\tfrac{z^n-z^{-n}}2,\tfrac{w^n+w^{-n}}2,\tfrac{w^n-w^{-n}}2, \]
\[\tfrac{(z^n+z^{-n})(w^n+w^{-n})}2,\tfrac{(z^n+z^{-n})(w^n-w^{-n})}2,
\tfrac{(z^n-z^{-n})(w^n+w^{-n})}2,\tfrac{(z^n-z^{-n})(w^n-w^{-n})}2,\]
where $n,m\geq1$.
In this case, colloquial notation  like
\[zw^{-1}-wz^{-1}\equiv\tfrac{(z-z^{-1})(w+w^{-1})}2-\tfrac{(z+z^{-1})(w-w^{-1})}2\]
is allowed.
There are natural $\mathfrak A[z]^\pm[w]^\pm$-module homomorphisms
$ \mathfrak A[z]^\pm[w]^\pm\rightarrow \mathfrak A \{ z, w\}
\rightarrow\mathfrak A\langle z\rangle\langle w\rangle$ respecting the grading.
\end{defin}
\begin{point}
The advantage of the terminology above is that it allows us to rewrite the definition of the
$\mathcal F$-square-root as the ``manifestly real'' expression
\[\sqrt[\mathcal F] T=\int\frac{(1+\frac{z+z^{-1}}2)T}{1+(z-2+z^{-1})T}\,\frac{|dz|}{2\pi}.\]
\end{point}
\begin{prop}If $T\in\spec_{\frac14-\overline{\mathbb R^-}}(\mathfrak A)$,
then $\frac{-T}{1-4T}$ is invertible and
$\sqrt[\mathcal F]{\frac{-T}{1-4T}}=\frac{-\sqrt[\mathcal F]T}{1-2\sqrt[\mathcal F]T}$.
Furthermore,
\begin{equation}\sqrt[\mathcal F] T(1-\sqrt[\mathcal F]T)=T.\label{e:sq2}\end{equation}
\begin{proof}
Substituting $z=-1$ into the resolvent term, we see that $(1-4T)^{-1}$ exists. The identity
$1+(z-2+z^{-1})\frac{-T}{1-4T} =(1-4T)^{-1}(1+((-z)-2+(-z)^{-1})T )$
shows that $\frac{-T}{1-4T}\in\spec_{\frac14-\overline{\mathbb R^-}}(\mathfrak A)$.
If \eqref{e:sq2} holds, then $(1-2\sqrt[\mathcal F]T)^2=1-4T$, and it is sufficient to prove that
$\sqrt[\mathcal F]{\frac{-T}{1-4T}}=\frac{\sqrt[\mathcal F]T-2T}{1-4T}$.
This, however, follows from the identity
\[\frac{(1+\frac{z+z^{-1}}2)\frac{-T}{1-4T}}{1+(z-2+z^{-1})\frac{-T}{1-4T}}
=\frac{\frac{(1+\frac{(-z)+(-z)^{-1}}2)T}{1+((-z)-2+(-z)^{-1})T}-2T}{1-4T}\]
integrated. So, what we have to show is the $\mathcal F$-square-root identity \eqref{e:sq2}.
Now, as
\begin{equation}\sqrt[\mathcal F] T(1-\sqrt[\mathcal F] T)-T=
\iint\!\frac{(1+\frac{z+z^{-1}}2)T}{1+(z-2+z^{-1})T}\left(1-
\frac{(1+\frac{w+w^{-1}}2)T}{1+(w-2+w^{-1})T}
\right)-T\frac{|dz|}{2\pi}\frac{|dw|}{2\pi},
\label{e16}\end{equation}
we should show that this latter term is $0$.
It would be natural to proceed along the steps
\begin{multline}
\frac{(1+\frac{z+z^{-1}}2)T}{1+(z-2+z^{-1})T}\left(1-
\frac{(1+\frac{w+w^{-1}}2)T}{1+(w-2+w^{-1})T}
\right)-T=\\
=\frac{\left(
\frac{z+z^{-1}}2-\frac{z+z^{-1}}2T-\frac{w+w^{-1}}2T
+T+\frac{z+z^{-1}}2\frac{w+w^{-1}}2T
\right)T(1-4T)}{(1+(z-2+z^{-1})T)(1+(w-2+w^{-1})T)}\\
\sim\frac{\left(
\left[\frac{\frac{z+z^{-1}}2+\frac{w+w^{-1}}2}2\right]
(1-2T)
+\left(1+\frac{z^{-1}w+zw^{-1}}2\right)T
\right)T(1-4T)}{(1+(z-2+z^{-1})T)(1+(w-2+w^{-1})T)}\\
=\left[\frac12\,\frac{z+w}{z-w}\right]
\frac{\left(\left(
\frac{z-z^{-1}}2-\frac{w-w^{-1}}2\right)(1-2T)
+(zw^{-1}-z^{-1}w)T
\right)T(1-4T)}{(1+(z-2+z^{-1})T)(1+(w-2+w^{-1})T)}\\
=\left[\frac12\,\frac{z+w}{z-w}\right]\left(
\frac{\frac{z-z^{-1}}2}{1+(z-2+z^{-1})T)}-\frac{\frac{w-w^{-1}}2}{1+(w-2+w^{-1})T)}
\right)T(1-4T)\sim 0,
\label{e17}\end{multline}
except we have to demonstrate that the use of $\frac12$ is superficial.

Another reordering operation can be defined according to
\[::z^nw^m::_{z,w}=\tfrac{z^{\max(|n|,|m|)}+z^{-\max(|n|,|m|)}}2
\tfrac{w^{\min(|n|,|m|)}+w^{-\min(|n|,|m|)}}2.\]
This applies to our standard $\mathfrak A[z]^\pm[w]^\pm$-modules.
In the context of  $\mathfrak A \langle z\rangle\langle w\rangle$,
it leaves only the $(++)$-graded parts and reorders them.
\begin{lemma} (a) For any Laurent series $a(z,w)\in \mathfrak A\langle z\rangle\langle w\rangle$, we have
\[\iint a(z,w)\frac{|dz|}{2\pi}\frac{|dw|}{2\pi}=\iint::a(z,w)::_{z,w}\frac{|dz|}{2\pi}\frac{|dw|}{2\pi}.\]

(b) If $::a_1(z,w)::_{z,w}=::a_2(z,w)::_{z,w}$ and
$b(z,w)\in\mathfrak A[z]^+[w]^+$ is symmetric, i. e. $b(z,w)=b(w,z)$, then
$::a_1(z,w) b(z,w)::_{z,w}=::a_2(z,w) b(z,w)::_{z,w}$.
\end{lemma}
For any symmetric Laurent series $b(z,w)\in \mathfrak A[z]^+[w]^+$, i. e. such that $b(z,w)=b(w,z)$,
we define formally
\[::\left[\frac{\frac{z+z^{-1}}2+\frac{w+w^{-1}}2}2\right]b(z,w)::_{z,w}=::\tfrac{z+z^{-1}}2b(z,w)::_{z,w}.\]

Let us consider a Laurent series $c(z,w)\in(\mathfrak A\{z,w\})^-$ such that it is antisymmetric, i. e.
$c(z,w)=-c(w,z)$. Then $c(z,w)$ is a formal linear combination of the basis elements
\[c_{n,0}=\tfrac{z^n-z^{-n}}2-\tfrac{w^n-w^{-n}}2\text{\quad and\quad}
c_{n,m}=\tfrac{(z^n-z^{-n})(w^m+w^{-m})}2-\tfrac{(w^n-w^{-n})(z^n+z^{-n})}2, \]
where $n,m\geq1$.
For such $c(z,w)$, we can formally  define $::\left[\frac12\,\frac{z+w}{z-w}\right]c(z,w)::_{z,w}$
according to our natural expectations. Again, we have to check that
the result is integral in terms of coefficients of $c(z,w)$.
We just give some samples in the table
\[\begin{array}{c||c|c|c|c}
::\left[\frac12\,\frac{z+w}{z-w}\right]c_{n,m}::_{z,w}&m=0&m=1&m=2&m=3\\
\hline
n=1&d^1_0&d^1_1+d^0_0&2d^1_0&-d^2_2+d^1_1+2d^2_0\\
n=2&d^2_0+d^1_1&2d^2_1+2d^1_0&d^2_2+2d^1_1+d^0_0&2d^2_1+2d^1_0\\
n=3&d^3_0+2d^2_1&2d^3_1+d^2_2+d^1_1+2d^2_0&2d^3_2+2d^2_1+2d^1_0&d^3_3+2d^2_2+2d^1_1+d^0_0
\end{array}\]
where $d^n_m=::z^nw^m::_{z,w}$.
In fact, at first sight, it looks more natural to choose $c(z,w)$ from
the antisymmetric elements of
$(\mathfrak A\langle z\rangle\langle w\rangle)^-=
\mathfrak A\langle z\rangle^+\langle w\rangle^-\oplus\mathfrak A\langle z\rangle^-\langle w\rangle^+$,
but it turns out that the coefficients
in $::\left[\frac12\,\frac{z+w}{z-w}\right]c(z,w)::_{z,w}$ would fail to be integral.
\begin{lemma} For any symmetric Laurent series $b(z,w)\in \mathfrak A[z]^+[w]^+$, we have
\[::\left[\frac{\frac{z+z^{-1}}2+\frac{w+w^{-1}}2}2\right]b(z,w)::_{z,w}=
::\left[\frac12\,\frac{z+w}{z-w}\right]\left(\tfrac{z-z^{-1}}2-\tfrac{w-w^{-1}}2\right) b(z,w)::_{z,w},\]
and
\[::\left(1+\frac{zw^{-1}+wz^{-1}}2\right)b(z,w)::_{z,w}=
::\left[\frac12\,\frac{z+w}{z-w}\right]\left(zw^{-1}-wz^{-1}\right) b(z,w)::_{z,w}.\notag\]
\end{lemma}
\begin{lemma} For any Laurent series $a(z)\in\mathfrak A\langle z\rangle^-$, we have
\[\iint::\left[\frac12\,\frac{z+w}{z-w}\right] (a(z)-a(w))::_{z,w}\frac{|dz|}{2\pi}\frac{|dw|}{2\pi}=0.\]
\end{lemma}
Now, it is easy to carry out the proof.
From \eqref{e16}, we should pass to the ``normal ordered'' form, after which the subsequent
manipulations leading to $0$ make sense.
\end{proof}
\end{prop}
\begin{prop}$P\in\spec_{\frac12+\mathrm i\overline{\mathbb R}}(\mathfrak A)$
if and only if
$P(1-P)\in\spec_{\frac14-\overline{\mathbb R^-}}(\mathfrak A)$.
Furthermore,
\[\sqrt[\mathcal F]{P(1-P)}=P+\idem P-2P\idem P.\]
Consequently,
\[\idem P =\frac{\sqrt[\mathcal F]{P(1-P)}-P}{1-2P}.\]
\begin{proof}
The decomposition
$1+(z-2+z^{-1})P(1-P)=(1-P+Pz)(1-P+Pz^{-1})$
implies the spectral statement.
The equality follows from the identity
\[\frac{\left(1+\frac{z+z^{-1}}2\right)P(1-P)}{1+(z-2+z^{-1})P(1-P)}
+\frac{\frac{z-z^{-1}}2P(1-P)(1-2P)}{1+(z-2+z^{-1})P(1-P)}=P+\frac{Pz}{1-P+Pz}(1-2P)\]
integrated.
\end{proof}
\end{prop}
Then we let $|P|_{\mathcal F}=\sqrt[\mathcal F]{P(1-P)}$.
Further statements can be proven parallel to the case with $\frac12$,
except the formulas are less customary.
E. g., the analogue of Proposition \ref{lem:sqex} is
\begin{prop} For $T\in\spec_{\frac14-\overline{\mathbb R^-}}(\mathfrak A)$, we have
\[\frac1{1+(z-2+z^{-1})T}=\frac1{1-2\sqrt[\mathcal F]T}\left(1
+\frac{-\sqrt[\mathcal F]T}{1-\sqrt[\mathcal F]T}(z+z^{-1})+
\left(\frac{-\sqrt[\mathcal F]T}{1-\sqrt[\mathcal F]T}\right)^2(z^2+z^{-2})+\ldots\right)\]
and
\[\frac1{1-\sqrt[\mathcal F]T}=\int\frac{1+z}{1+(z-2+z^{-1})T}\frac{|dz|}{2\pi}.\]
\end{prop}

\end{document}